\documentclass[letterpaper,11pt]{amsart}

\usepackage{latexsym,array,delarray,amsthm,amssymb,epsfig}%eufrak

\theoremstyle{plain}
\newtheorem{thm}{Theorem}[section]
\newtheorem{lemma}[thm]{Lemma}
\newtheorem{prop}[thm]{Proposition}
\newtheorem{cor}[thm]{Corollary}

\theoremstyle{definition}
\newtheorem{df}[thm]{Definition}
\newtheorem{ex}[thm]{Example}

\newcommand{\zz}{\mathbb{Z}}
\newcommand{\nn}{\mathbb{N}}
\newcommand{\pp}{\mathbb{P}}

\newcommand{\qq}{\mathbb{Q}}
\newcommand{\rr}{\mathbb{R}}

\newcommand{\kk}{\mathbb{K}}
\newcommand{\aaa}{\mathbf{a}}
\newcommand{\bbb}{\mathbf{b}}

\newcommand{\ca}{\mathcal{A}}

\begin{document}

\title{Toric Fiber Products}

\author{Seth Sullivant}
\address{Department of Mathematics, Harvard University, Cambridge, MA 02138}
\email{seths@math.harvard.edu}

\maketitle

\begin{abstract}
We introduce and study the toric fiber product of two ideals in polynomial rings that are homogeneous with respect to the same multigrading.  Under the assumption that the set of degrees of the variables form a linearly independent set, we can explicitly describe generating sets and Gr\"obner bases for these ideals.  This allows us to unify and generalize some results in algebraic statistics.
\end{abstract}

%%%%%%%%%%%%%%%%%%%%%%%%%%%%%%%%%%%%%%%%
%%%%%%%%%%%%%%%%%%%%%%%%%%%%%%%%%%%%%%%%
%%%%%%%%%%%%%%%%%%%%%%%%%%%%%%%%%%%%%%%%
%%%%%%%%%%%%%%%%%%%%%%%%%%%%%%%%%%%%%%%%
%%%%%%%%%%%%%%%%%%%%%%%%%%%%%%%%%%%%%%%%

\section{Introduction}\label{sec:intro}

A common problem in algebraic statistics is to convert the parametric representation of a statistical model into the implicit representation in terms of finding the defining prime ideal of the model.  This is a special case of the implicitization problem that arises frequently in computational algebraic geometry.  In algebraic statistics, we are usually presented with a family of statistical models and we would like to find a theorem which gives a complete description of all the ideals for all the statistical models in this family.  A useful approach has been to try to find decomposition rules for the models and the resulting ideals, and subsequently reduce the problem to finding the defining prime ideals in a few special cases which can then be handled theoretically or using a computer algebra system.  This approach has played a role in attacking the problem of determining phylogenetic invariants for various tree-based models of evolution \cite{Allman2004, Sturmfels2005} and for studying Markov bases of hierarchical models \cite{Dobra2004, Hosten2002}.

In this paper, we introduce the toric fiber product, an operation that takes two homogeneous ideals with compatible multigradings and produces a new homogeneous ideal.  This operation generalizes  the Segre product of two schemes as well as the gluing operations for toric ideals that appear in \cite{Dobra2004, Hosten2002, Sturmfels2005}.  When the underlying grading group has special structure, we are able to explicitly compute generating sets and special Gr\"obner bases for the toric fiber product from generating sets and Gr\"obner bases of the component ideals.

Given a positive integer $n$,  let $[n] = \{1,2,\ldots, n\}$ denote the set of the first $n$ positive integers.  Let $r > 0$ be a positive integer and $s, t \in \zz_{> 0}^r$ be two vectors of positive integers.
Let $$\kk[x] = \kk[x^i_j \, \, | \, \, i \in [r], j \in [s_i] ]$$ and $$ \kk[y] = \kk[y^i_k \, \, | \, \, i \in [r], k \in [ t_i] ]$$ be multigraded polynomial rings subject to the multigrading 
$$\deg(x^i_j) = \deg(y^i_k) = \aaa^i \in \zz^d.$$  
We assume throughout that there exists a vector $\omega \in \qq^d$ such that $\omega^T\aaa^i = 1$ for all $i$. This implies that ideals in $\kk[x]$ or $\kk[y]$ that are homogeneous with respect to the multigrading are homogeneous in the usual sense.  Denote by $\ca = \{\aaa^1, \ldots, \aaa^r \}$ and let $\nn\ca$ be the affine semigroup generated  by $\ca$. 

If $I$ and $J$ are homogeneous ideals (with respect to the multigrading) in $\kk[x]$ and $\kk[y]$, respectively, the quotient rings $R = \kk[x]/I$ and $S = \kk[y]/J$ are also multigraded rings. Let 
$$\kk[z] = \kk[z^i_{jk} \, \, | \, \, i \in [r], j \in [s_i], k \in [t_i]]$$ 
and let $\phi_{I,J} : \kk[z] \rightarrow R \otimes_\kk S$ be the ring homomorphism such that $z^i_{jk} \mapsto x^i_j \otimes y^i_k$. 

\begin{df}
The \emph{toric fiber product} of $I$ and $J$, denoted $I \times_\ca J$ is the kernel of $\phi_{I,J}$:
$$ I \times_\ca J = \ker(\phi_{I,J}).$$
\end{df}

Two fundamental examples, illustrating the coarsest and finest possible multigradings, are Segre products and sums of monomial ideals.

\begin{ex}
Suppose that $r = 1$ and $I = J = 0$.  Then $\phi_{I,J}$ is the ring homomorphism
$$\phi_{I,J}: \kk[z_{jk} \, \, | \, \, j \in [s], k \in [t] ] \rightarrow \kk[x_j, y_k \, \, | \, \, j \in [s], k \in [t] ]$$
$$ z_{jk} \mapsto x_j y_k.$$
The toric fiber product is 
$$I \times_\ca J = \ker(\phi_{I,J}) = \langle  z_{j_1 k_1}z_{j_2 k_2} - z_{j_1 k_2}z_{j_2 k_1} \, \, | \, \, j_1, j_2 \in [s], k_1, k_2 \in [t] \rangle;$$
that is, $I \times_\ca J$ is the ideal of $2 \times 2 $ minors of a generic matrix.
\end{ex}

\begin{ex}\label{ex:mono}
Suppose that $r>0$ and $s, t = (1,\ldots, 1)$ are the all ones vector.   Suppose that $\deg(x^i) = \deg(y^i) = \deg(z^i) = e_i$, the $i$th standard unit vector.
If $I \in \kk[x]$ and $J \in \kk[y]$ are homogeneous with respect to this multigrading, they must both be monomial ideals.  Then the toric fiber product is simply
$$I \times_\ca J = I(z) + J(z)$$
where $I(z)$ denotes the ideal $I$ with $z$ variables substituted for the $x$ variables (and similarly for $J(z)$). 
\end{ex}

Our main interest in toric fiber products is when $I$ and $J$ are the prime ideals of unirational varieties.  Then the ideal $I \times_\ca J$ is also prime and defines a unirational variety.  In practice, we are often presented with the parametrization of a unirational variety and we are interested in finding its defining ideal.  One useful tool is to find a nice grading such that the ideal is, in fact, a toric fiber product.  This grading is usually considerably coarser than the finest grading associated to the ideal.  If $\ca$ is a linearly independent set,  we can determine generators and Gr\"obner bases of the toric fiber product $I \times_\ca J$ explicitly from generators and Gr\"obner bases of $I$ and $J$.

\begin{ex}\label{ex:dets}
Let $\phi$ be the ring homomorphism
$$\phi: \kk[q_{i_1\bullet i_3 \bullet i_5} \, \, | \, \, i_1,i_3,i_5 \in [3] ] \rightarrow 
\kk[a_{i_1i_2}, b_{i_2i_3}, c_{i_3i_4}, d_{i_4i_5} \, \, | \, \, i_1,i_3,i_5 \in [3], i_2,i_4 \in [2] ]$$
$$q_{i_1\bullet i_3\bullet i_5} \mapsto \sum_{i_2 = 1}^2  \sum_{i_4 = 1}^2  a_{i_1i_2}b_{i_2i_3} c_{i_3i_4} d_{i_4i_5},$$
Each of the polynomials appearing in the parametrization is homogeneous with respect to the grading $\deg(p_{i_1 \bullet i_3 \bullet i_5}) = e_{i_3}$, the $i_3$-th standard unit vector.
The polynomials appearing in the parametrization can be written in factored form as 
$$q_{i_1\bullet i_3 \bullet i_5} \mapsto \left( \sum_{i_2 =1}^2  a_{i_1i_2}b_{i_2i_3} \right) \left( \sum_{i_4 = 1}^2 c_{i_3i_4} d_{i_4i_5} \right).$$
The ideal $K = \ker( \phi )$ is a toric fiber product $K = I \times_\ca I$.  The underlying ring
$R =  \kk[q_{i_1\bullet i_3} \, \, | \, \, i_1,i_3 \in [3] ]$ has the grading $\deg(q_{i_1i_3}) = e_{i_3}$, the $k$th standard unit vector.  The ideal $I$ is the kernel of the ring homomorphism
$$\hat{\phi} :  R \rightarrow \kk[a_{i_1i_2}, b_{i_2i_3} \,\, | \, \, i_1,i_3 \in [3], i_2 \in [2] ]$$
$$ q_{i_1 \bullet i_3} \mapsto \sum_{i_2 = 1}^2 a_{i_1i_2} b_{i_2i_3}.$$
Thus, $I$ is the principal ideal generated by the  determinant of the matrix:
$$\begin{pmatrix}
q_{1\bullet 1} & q_{1\bullet 2} & q_{1\bullet 3} \\
q_{2\bullet 1} & q_{2\bullet 2} & q_{2\bullet 3} \\
q_{3\bullet 1} & q_{3\bullet 2} & q_{3\bullet 3}
\end{pmatrix}.$$

Using the machinery in Section \ref{sec:mono}, one can show that the ideal $K = I \times_\ca I$ is generated by determinants of \emph{flattenings} and \emph{slices} of the $3$-dimensional tensor $(q_{i_1\bullet i_3 \bullet i_5})$.  In particular, $K$ is generated by the $3 \times 3 $ minors of the matrices 
$$\begin{pmatrix}
q_{1\bullet 1\bullet 1} & q_{1\bullet 1\bullet 2} &  q_{1\bullet 1\bullet 3} & q_{1\bullet 2\bullet 1} & q_{1\bullet 2\bullet 2} &  q_{1\bullet 2\bullet 3} & q_{1\bullet 3\bullet 1} & q_{1\bullet 3\bullet 2} &  q_{1\bullet 3\bullet 3} \\
q_{2\bullet 1\bullet 1} & q_{2\bullet 1\bullet 2} &  q_{2\bullet 1\bullet 3} & q_{2\bullet 2\bullet 1} & q_{2\bullet 2\bullet 2} &  q_{2\bullet 2\bullet 3} & q_{2\bullet 3\bullet 1} & q_{2\bullet 3\bullet 2} &  q_{2\bullet 3\bullet 3} \\
q_{3\bullet 1\bullet 1} & q_{3\bullet 1\bullet 2} &  q_{3\bullet 1\bullet 3} & q_{3\bullet 2\bullet 1} & q_{3\bullet 2\bullet 2} &  q_{3\bullet 2\bullet 3} & q_{3\bullet 3\bullet 1} & q_{3\bullet 3\bullet 2} &  q_{3\bullet 3\bullet 3} \\
\end{pmatrix}
$$
$$
\begin{pmatrix}
q_{1\bullet 1\bullet 1} & q_{1\bullet 2\bullet 1} &  q_{1\bullet 3\bullet 1} & q_{2\bullet 1\bullet 1} & q_{2\bullet 2\bullet 1} &  q_{2\bullet 3\bullet 1} & q_{3\bullet 1\bullet 1} & q_{3\bullet 2\bullet 1} &  q_{3\bullet 3\bullet 1} \\
q_{1\bullet 1\bullet 2} & q_{1\bullet 2\bullet 2} &  q_{1\bullet 3\bullet 2} & q_{2\bullet 1\bullet 2} & q_{2\bullet 2\bullet 2} &  q_{2\bullet 3\bullet 2} & q_{3\bullet 1\bullet 2} & q_{3\bullet 2\bullet 2} &  q_{3\bullet 3\bullet 2} \\
q_{1\bullet 1\bullet 3} & q_{1\bullet 2\bullet 3} &  q_{1\bullet 3\bullet 3} & q_{2\bullet 1\bullet 3} & q_{2\bullet 2\bullet 3} &  q_{2\bullet 3\bullet 3} & q_{3\bullet 1\bullet 3} & q_{3\bullet 2\bullet 3} &  q_{3\bullet 3\bullet 3} \\
\end{pmatrix}
$$
together with the $2 \times 2$ minors of the matrices
$$
\begin{pmatrix}
q_{1\bullet 1\bullet 1} & q_{1\bullet 1\bullet 2} & q_{1\bullet 1\bullet 3} \\
q_{2\bullet 1\bullet 1} & q_{2\bullet 1\bullet 2} & q_{2\bullet 1\bullet 3} \\
q_{3\bullet 1\bullet 1} & q_{3\bullet 1\bullet 2} & q_{3\bullet 1\bullet 3} \\
\end{pmatrix}
\quad
\begin{pmatrix}
q_{1\bullet 2\bullet 1} & q_{1\bullet 2\bullet 2} & q_{1\bullet 2\bullet 3} \\
q_{2\bullet 2\bullet 1} & q_{2\bullet 2\bullet 2} & q_{2\bullet 2\bullet 3} \\
q_{3\bullet 2\bullet 1} & q_{3\bullet 2\bullet 2} & q_{3\bullet 2\bullet 3} \\
\end{pmatrix}
$$
$$
\begin{pmatrix}
q_{1\bullet 3\bullet 1} & q_{1\bullet 3\bullet 2} & q_{1\bullet 3\bullet 3} \\
q_{2\bullet 3\bullet 1} & q_{2\bullet 3\bullet 2} & q_{2\bullet 3\bullet 3} \\
q_{3\bullet 3\bullet 1} & q_{3\bullet 3\bullet 2} & q_{3\bullet 3\bullet 3} \\
\end{pmatrix}.
$$
Furthermore, this collection of $2 \times 2$ and $3 \times 3$ minors form a Gr\"obner basis for $K$. 
This example is a special case of Corollary \ref{cor:chain}.
\qed
\end{ex}

The main focus of this paper is on the special case of toric fiber products where $\ca$ is a linearly independent set.  As we will see, this played a significant role in Example \ref{ex:dets}.  
It is probably impossible to recover explicitly the generating set of $K = I \times_\ca J$ from the ideals $I$ and $J$, if $\ca$ is not linearly independent.  Indeed, in the next section, we will see an example where $I$ is generated by quadrics and $J = \langle 0 \rangle$ but $I \times_\ca J$ requires minimal generators of arbitrarily large degree.

The outline for this paper is as follows.  In the next section, we use a contraction of an ideal under a monomial homomorphism to determine generating sets and Gr\"obner bases for the toric fiber products $I \times_\ca J$, when $\ca$ is linearly independent.    In Section \ref{sec:apps}, we consider applications of the main result.  This includes proofs of some known results about the defining ideals of products of projective schemes, and their Gr\"obner bases.  We also illustrate how the toric fiber product arises in algebraic statistics.  This allows us to unify results about reducible hierarchical models and group-based models on phylogenetic trees.

%%%%%%%%%%%%%%%%%%%%%%%%%%%%%%%%%%%%%%%%
%%%%%%%%%%%%%%%%%%%%%%%%%%%%%%%%%%%%%%%%
%%%%%%%%%%%%%%%%%%%%%%%%%%%%%%%%%%%%%%%%
%%%%%%%%%%%%%%%%%%%%%%%%%%%%%%%%%%%%%%%%

\section{Contractions Under Monomial Homomorphisms}\label{sec:mono}

Let $B \in \zz^{d\times n}$ be a $d \times n$ integral matrix.  Consider the ring homomorphism
$$\phi_B : \kk[z_1, \ldots, z_n] \rightarrow  \kk[t_1, \ldots, t_d]$$
$$ z_j \mapsto  \prod_{i = 1}^d  t_i^{b_{ij}}.$$
We call $\phi_B$ a monomial homomorphism.  The ideal $I_B = \ker(\phi_B)$ is called a \emph{toric ideal}.  Sturmfels' book \cite{Sturmfels1996} is a standard reference for background on toric ideals.  In this section, we consider the contractions $\phi_B^{-1}(I)$ of arbitrary ideals $I$ in $\kk[t]$, and apply these results to toric fiber products.  The main idea here is to compare the initial ideals of $I$ to the initial ideals of $\phi_B^{-1}(I)$.

Let $\omega \in \zz_{\geq 0}^d$ be a vector of weights.  The vector $\omega$ induces a partial order $\prec_\omega$ on the set of monomials in $\kk[t]$ by declaring that $t^\aaa \prec_\omega t^\bbb$ if $\omega^T \aaa < \omega^T \bbb$.    The partial order $\prec_\omega$ is called a \emph{weight order} on $\kk[t]$.  Note that it need not be a term order on $\kk[t]$.  

Given a polynomial $f \in \kk[t]$, the initial form ${\rm in}_\omega(f)$ is the sum of all terms of $f$ that have the highest weight with respect to the partial order $\prec_\omega$.  If $I$ is an ideal of $\kk[t]$, the initial ideal ${\rm in}_\omega(I)$ is the ideal
$${\rm in}_\omega (I) = \left<  {\rm in}_\omega(f) \, \, | \, \, f \in I \right>.$$
Our main use for for weight orders comes from the following useful fact.

\begin{prop}\cite[Prop. 1.11]{Sturmfels1996}
For any term order $\prec$ and any ideal $I \subset \kk[t]$ there exists a vector $\omega \in \zz^d_{\geq 0}$ such that ${\rm in}_\omega (I) = {\rm in}_\prec (I)$.
\end{prop}

We say that a finite collection of polynomials $G \subset I$ is a Gr\"obner basis of $I$ with respect to the weight order $\omega$ if $\langle {\rm in}_\omega (g) \, | \, g \in G \rangle = {\rm in}_\omega(I)$ and this ideal is a monomial ideal.  A collection of polynomials such that $\langle {\rm in}_\omega (g) \, | \, g \in G \rangle = {\rm in}_\omega(I)$ is called a pseudo-Gr\"obner basis.  Any Gr\"obner basis or pseudo-Gr\"obner basis of $I$ generates $I$.

Every weight order $\prec_\omega$ on $\kk[t]$ determines a weight order $\prec_{\phi_B^* \omega}$ on $\kk[z]$ via the pullback of $\omega$ through $\phi$.  That is, $\phi_B^*\omega = \omega^T B$.  Note that this construction has two important properties.  First, if $z^\aaa$ is a monomial then its weight with respect to $\phi_B^* \omega$ equals the weight of $\phi_B(z^\aaa)$ with respect to $\omega$.  Second, if $f \in \ker(\phi_B)$ is a binomial then ${\rm in}_{\phi^*_B \omega}(f) = f$, which implies that ${\rm in}_{\phi^*_B \omega}(I_B) = I_B$ for all $\omega \in \rr^d$.

\begin{lemma}\label{lem:initials}
Let $I$ be an ideal in $\kk[t]$.  Then
$$ {\rm in}_{\phi^*_B \omega}( \phi^{-1}_B(I) )   \subseteq  \phi^{-1}_B( {\rm in}_{\omega} (I)).$$
\end{lemma}

\begin{proof}
Let $f \in \phi^{-1}_B(I)$.  We must show that $\phi_B({\rm in}_{\phi^*_B \omega}(f)) \in {\rm in}_\omega(I)$.  Without loss of generality, we may assume that $f$ is reduced with respect to any Gr\"obner basis of the toric ideal $I_B$.  This is because $I_B \subset {\rm in}_{\phi^*_B \omega}(\phi^{-1}(I) ) $ for any $I$ and $\phi_B(I_B) = 0$.  Since $f$ is reduced with respect to $I_B$, this means there is at most one term of $f$ in each $B$ graded degree.  In particular, there can be no cancellation amongst the terms of $\phi_B(f)$.  Since each pair of monomial of $f$ and image monomial of $\phi_B(f)$ have the same weight with respect to $\phi^*_B \omega$ and $\omega$ respectively, we deduce that $\phi_B({\rm in}_{\phi^*_B \omega}(f)) \in {\rm in}_\omega(I)$ which completes the proof.
\end{proof}

\begin{lemma}\label{lem:monosplit}
Let $M = \left< m_1, \ldots, m_r \right> \subset \kk[t]$ be a monomial ideal.  Then 
$$\phi^{-1}_B(M) = \phi^{-1}_B( \left< m_1 \right> ) + \cdots + \phi^{-1}_B( \left< m_r \right> ).$$
Furthermore $\phi^{-1}_B(M) = M' + I_B$ where $M'$ is a monomial ideal.
\end{lemma}

\begin{proof}
We use the same argument as the proof of Lemma \ref{lem:initials}.  In particular, suppose $f \in \phi^{-1}_B(M)$.  We may suppose that $f$ is reduced with respect to any Gr\"obner basis of $I_B$.  Then each monomial in $f$ maps to a monomial in $\kk[t]$, and $f \in \phi^{-1}_B(M)$ if and only if each monomial of $\phi_B(f)$ belongs to $M$.  This means that if $n$ is a monomial of $f$, $\phi_B(n) \in \left< m_i \right>$ for some $i$.  Furthermore, this shows that each monomial of $f$ belongs to $\phi^{-1}_B(M)$.   
\end{proof}

Lemmas \ref{lem:initials} and \ref{lem:monosplit} suggest a strategy for determining the ideals $\phi^{-1}_B(I)$.  First, we compute an initial ideal of $I$.  Then we determine $\phi^{-1}_B ({\rm in}_\omega (I))$ using combinatorial arguments. Then, if we are lucky, we find a collection of polynomials $G \subset \phi^{-1}_B(I)$ such that $\left< {\rm in}_{\phi^*_B \omega}(G) \right> = {\rm in}_{\phi^*_B \omega}(\phi^{-1}_B(I))$.  Then we can conclude that $G$ is a pseudo-Gr\"obner basis for $\phi^{-1}_B(I)$ with respect to the weight order $\prec_{\phi^*_B \omega}$.

In general, this strategy is not possible to implement, because either
$ {\rm in}_{\phi^*_B \omega}( \phi^{-1}_B(I) )   \neq  \phi^{-1}_B( {\rm in}_{\omega} (I))$ or there is no combinatorial description of $\phi^{-1}_B(M)$ where $M$ is a monomial ideal, or both.  However, in the special case that arises when taking a toric fiber product, we will see that there is a simple answer to both problems.

To give our main algebraic result concerning the Gr\"obner bases and generating sets of the toric fiber products $I \times_\ca J$, we first need to show that the toric fiber product is the contraction of a monomial homomorphism $\phi^{-1}_B(I+J)$, for suitable $B$.  To this end, we will derive an alternate description of the toric fiber product which fits into this framework.

With the setup from Section \ref{sec:intro}, consider the monomial homomorphism 
$$\phi_B: \kk[z] \rightarrow \kk[x,y] := \kk[x^i_j, y^i_k \, \, | \, \, i \in [r], j \in [s_i], k \in [t_i ] ]$$
\begin{equation}\label{eqn:B}
 z^i_{jk} \mapsto x^i_j y^i_k.
 \end{equation}

For the remainder of this section $B$ denotes the matrix arising from a toric fiber product according to Equation \ref{eqn:B}.  This matrix only depends on $r$, $s$, and $t$.

\begin{prop} 
$$I \times_\ca J = \phi^{-1}_B(I + J).$$
\end{prop}
Note that  the ideal $I +J$ is considered as an ideal in $\kk[x,y]$ after taking the extensions of $I$ and $J$.
\begin{proof}
Note that $\kk[x,y] \cong \kk[x] \otimes_\kk \kk[y]$.  Hence $\kk[x,y] / (I + J) \cong \kk[x]/I \otimes_\kk \kk[y]/J$.  Given any ring homomorphism $\phi: R \rightarrow S$ and induced ring homomorphism $\hat{\phi}: R \rightarrow S/I$, we have $\ker(\hat{\phi}) = \phi^{-1}(I)$, and this completes the proof.
\end{proof}

\begin{lemma}\label{lem:monopull}
Let $m = x^{i_1}_{j_1} x^{i_2}_{j_2} \cdots x^{i_d}_{j_d}$ be a monomial in $\kk[x,y]$.  Then 
$$\phi^{-1}_B(\left<m\right>) = \left<  z^{i_1}_{j_1k_1} z^{i_2}_{j_2 k_2} \cdots z^{i_d}_{j_dk_d} \, \, | \, \, k_1 \in [t_{i_1}], \ldots k_d \in [t_{i_d}] \right>  + I_B.$$
Similarly, if $n =   y^{i_1}_{k_1} y^{i_2}_{k_2} \cdots y^{i_d}_{k_d}$ is a monomial in $\kk[x,y]$, then
$$\phi^{-1}_B(\left< n \right> ) = \left<   z^{i_1}_{j_1k_1} z^{i_2}_{j_2 k_2} \cdots z^{i_d}_{j_dk_d} \, \, | \, \, 
j_1 \in [s_{i_1}], \ldots j_d \in [s_{i_d}] \right>  + I_B.$$
\end{lemma}

\begin{proof}
It suffices to handle the first case.  By the second part of  Lemma \ref{lem:monosplit} it suffices to determine the monomials that belong to $\phi^{-1}_B(\left<m\right>)$.  Denote by $M'$ the monomial ideal on the right hand side of the equation.  Given a monomial $m' =  \prod z^{i_l}_{j_lk_l} \in \kk[z]$, its image is the monomial $\prod x^{i_l}_{j_l} y^{i_l}_{k_l}$ which belongs to $\left< m \right>$ if and only if $\prod x^{i_l}_{j_l}$ belongs to $\left< m \right>$.  But this implies that these exists a monomial in $M'$ dividing $m'$.
\end{proof}

\begin{prop}\label{prop:quad}
Let $B$ be the matrix representing the monomial homomorphism arising from the toric fiber product.  Then
$$I_B = \left< \underline{z^i_{j_1k_2} z^i_{j_2k_1}} - z^i_{j_1k_1} z^i_{j_2k_2} \, \, | \, \, 
1 \leq i \leq r, 1 \leq j_1 < j_2 \leq s_i, 1 \leq k_1 < k_2 \leq t_i   \right> $$
and these quadrics are a Gr\"obner basis for $I_B$ with respect to any term order that selects the underlined terms as leading terms.
\end{prop}

Denote by ${\rm Quad}_B$ the set of quadrics described in Proposition \ref{prop:quad}.

\begin{proof}
First of all, there exist term orders which select the underlined terms as leading terms.   Indeed,
let $\prec$ be the lexicographic term order such that
$z^{i_1}_{j_1k_1} \prec z^{i_2}_{j_2k_2}$ if $i_1 < i_2$ or $i_1 = i_2$ and $j_1 < j_2$ or $i_1 = i_2$ and $j_1 = j_2$ and $k_1 > k_2$.  In particular, $\prec$ selects the underlined terms of the quadrics in ${\rm Quad}_B$ as leading terms.

Since $I_B$ is a toric ideal, it suffices to show that if $f$ is any binomial in $I_B$, there exists a quadric  $g \in {\rm Quad}_B$ such that ${\rm in}_\prec (g) | {\rm in}_\prec (f)$.  To each binomial 
$$f = z^{i_1}_{j_1k_1} \cdots z^{i_d}_{j_d k_d} - z^{i'_1}_{j'_1 k'_1} \cdots z^{i'_d}_{j'_d k'_d}$$
we associate the tableaux of indices
$$f = \left[ \begin{array}{ccc}
i_1 & j_1 & k_1 \\
\vdots & \vdots & \vdots \\
i_d & j_d & k_d
\end{array} \right]
-
\left[ \begin{array}{ccc}
i'_1 & j'_1 & k'_1 \\
\vdots & \vdots & \vdots \\
i'_d & j'_d & k'_d
\end{array} \right].
$$
Note that two individual tableau represent the same monomial if and only if one can be obtained from the other by swapping rows.  A binomial $f$ belongs to $I_B$ if and only if the image of each monomial under $\phi_B$ is the same.  In tableau notation, this can be expressed as
$$ 
\left[ \begin{array}{cc}
i_1 & j_1 \\
\vdots & \vdots \\
i_d & j_d \end{array} \right]
\left[ \begin{array}{cc}
i_1 & k_1 \\
\vdots & \vdots \\
i_d & k_d \end{array} \right]
\quad = 
\quad
\left[ \begin{array}{cc}
i'_1 & j'_1 \\
\vdots & \vdots \\
i'_d & j'_d \end{array} \right]
\left[ \begin{array}{cc}
i'_1 & k'_1 \\
\vdots & \vdots \\
i'_d & k'_d \end{array} \right].
$$
In each expression, the first tableau denotes the indices of the $x$ variables and the second tableau denotes the indices of the $y$ variables.  Thus, after rearranging the rows of the second tableau of $f$ we may write
$$f = \left[ \begin{array}{ccc}
i_1 & j_1 & k_1 \\
\vdots & \vdots & \vdots \\
i_d & j_d & k_d
\end{array} \right]
-
\left[ \begin{array}{ccc}
i_1 & j_1 & k'_1 \\
\vdots & \vdots & \vdots \\
i_d & j_d & k'_d
\end{array} \right]
$$
Furthermore, we will assume that the rows of each tableau are ordered lexicographically with $111 \prec 112 \prec \cdots$.  A monomial is reduced with respect to the set ${\rm Quad}_B$ if and only if there is no sub-tableau of the form
$$
\left[ \begin{array}{ccc}
i & j_1 & k_2 \\
i & j_2  & k_1
\end{array} \right]
$$
where $j_1 < j_2$ and $k_1 < k_2$.  This implies that in a monomial that has been reduced with respect to ${\rm Quad}_B$, for each fixed value of $i$, both the $j$'s and $k$'s are nondecreasing down their column in the subtableau containing all the rows with $i_t = i$.  However, among all the tableaux that have the same image under $\phi_B$, there is only one that has this property and it is minimal with respect to $\prec$.  Thus, if $f \in I_B$ its leading term must be divisible by some leading term in ${\rm Quad}_B$.
\end{proof}

For the rest of this section, we will assume that $\ca$ is a linearly independent set.  Let $f \in I$ be a homogeneous polynomial (with respect to the multigrading by $\nn \ca$).  If $f$ has degree $d$ and $v$ terms we can write
$$ f = \sum_{u = 1}^v c_u x^{i^u_1}_{j^u_1} x^{i^u_2}_{j^u_2} \cdots x^{i^u_d}_{j^u_d},$$
where each $c_u \in \kk$ is a coefficient.  However, the fact that $f$ is homogeneous and $\ca$ is linearly independent guarantees that each multiset of upper indices $M_u = \{ i^u_1, \ldots, i^u_d \}$ is independent of $u$.  That is, for all $u, u'$, $M_u = M_{u'}$.  Thus after possibly rearranging the indeterminates appearing in each monomial we can always write
$$f = \sum_{u =1}^v c_u x^{i_1}_{j^u_1} x^{i_2}_{j^u_2} \cdots x^{i_d}_{j^u_d}.$$
Now let $k = (k_1, \ldots, k_d)$ with $k_1 \in [t_{i_1}], k_2, \in [t_{i_2}], \ldots, k_d \in [t_{i_d}]$ and consider the polynomial $f_k \in \kk[z]$ defined by
$$f_k = \sum_{u = 1}^v  c_u  z^{i_1}_{j^u_1 k_1} z^{i_2}_{j^u_2 k_2 }  \cdots  z^{i_d}_{j^u_d k_d}.$$

Note that since $f \in I$, the new homogeneous polynomial $f_k \in I \times_\ca J$ for all $k$.  This follows because 
$$\phi_B(f_k)  =  y^{i_1}_{k_1} y^{i_2}_{k_2} \cdots y^{i_d}_{k_d} \cdot f \in I.$$

\begin{df}
Let $\ca$ be linearly independent and let $F \subset I$ be a collection of homogeneous polynomials.  To each $f \in F$ we associate the set 
$T_f =  \prod_{l =1}^d [t_{i_l}]$ of indices.  Denote by
$${\rm Lift}(F) \, = \, \left\{ f_k \, \, | \, \, f \in F,  k \in T_f \right\}$$
which we call the lifting of $F$ to $I \times_\ca J$.  If $G \subset J$ is a collection of homogeneous polynomials we define ${\rm Lift}(G)$ in the analogous way.
\end{df}

\begin{thm}\label{thm:main}
Suppose that $\ca$ is linearly independent.  Let $F \subset I$ be a homogeneous Gr\"obner basis for $I$ with respect to the weight vector $\omega_1$ and let $G \subset J$ be a homogeneous Gr\"obner basis for $J$ with respect to the weight vector $\omega_2$.  Then 
$${\rm Lift}(F) \cup {\rm Lift}(G) \cup {\rm Quad}_B$$
is a pseudo-Gr\"obner basis for $I \times_\ca J$ with respect to the weight order $\phi^*_B(\omega_1, \omega_2)$.
\end{thm}

\begin{proof}
First of all, note that since the generators of $I$ and $J$ are in disjoint sets of variables, $F \cup G$ is a Gr\"obner basis for $I +J$ with respect to the weight order $(\omega_1,\omega_2)$.  Let $M = {\rm in}_{(\omega_1,\omega_2)}(I+J)$ be the initial ideal of $I +J$ with respect to $(\omega_1,\omega_2)$.    Since $M$ is a monomial ideal and each minimal generator of $I + J$ belongs to $\kk[x]$ or $\kk[y]$, we can use Lemmas \ref{lem:monosplit} and \ref{lem:monopull} to compute $\phi^{-1}_B(M)$.  However, each of the monomials appearing in $\phi^{-1}_B(M)$ appears as the $\phi^*_B(\omega_1, \omega_2)$ leading term of some polynomial in ${\rm Lift}(F)$ or ${\rm Lift}(G)$.   Additionally, ${\rm Quad}_B$ generates $I_B$.  So we deduce that the initial forms in ${\rm Lift}(F) \cup {\rm Lift}(G) \cup {\rm Quad}_B \subset I \times_\ca J$ generate $\phi^{-1}_B(M)$.  Thus, by Lemma \ref{lem:initials}, we deduce that $\phi^{-1}_B(M) = {\rm in}_{\phi^*_B(\omega_1,\omega_2)}(I \times_\ca J)$ and this completes the proof.
\end{proof}

\begin{thm}\label{thm:maingro}
With the same assumptions as Theorem \ref{thm:main}, let $\omega$ be a weight vector such that ${\rm Quad}_B$ is a Gr\"obner basis for $I_B$.  Then
$${\rm Lift}(F) \cup {\rm Lift}(G) \cup {\rm Quad}_B$$
is a Gr\"obner basis for $I \times_\ca J$ with respect to the weight order $\phi^*_B(\omega_1, \omega_2) + \epsilon \omega$ for sufficiently small $\epsilon > 0$. 
\end{thm}

\begin{proof}
If we choose $\epsilon$ very small, we will have ${\rm in}_{\phi^*_B(\omega_1,\omega_2)}(f_k) = {\rm in}_{\phi^*_B(\omega_1,\omega_2) + \epsilon \omega}(f_k) $ for all $f_k \in {\rm Lift}(F) \cup {\rm Lift}(G)$.  This implies that ${\rm in}_{\phi^*_B(\omega_1, \omega_2) + \epsilon \omega}(I \times_\ca J) = {\rm in}_\omega(  {\rm in}_{\phi^*_B(\omega_1, \omega_2)}( I \times_\ca J))$ since we only need to determine the initial terms of $\left< {\rm Quad}_B \right>$.  But because $\omega$ induces a term order that make ${\rm Quad}_B$ into a Gr\"obner basis for $I_B$ we are done.
\end{proof}

\begin{cor}
Let $\ca$ be a linearly independent set.  Let $F$ be a homogeneous generating set for $I$ and $G$ be a homogeneous generating set for $J$.  Then
$${\rm Lift}(F) \cup {\rm Lift}(G) \cup {\rm Quad}_B$$
is a generating set for $I \times_\ca J$.
\end{cor}

\begin{proof}
If $F$ and $F'$ generate $I$, ${\rm Lift}(F)$ and ${\rm Lift} (F')$ generate the same ideal.  This holds, in particular, if $F'$ is a Gr\"obner basis for $I$.  A similar statement holds of $G$ as well.  Thus
\begin{equation}\label{eqn:gens}
\left< {\rm Lift}(F) \cup {\rm Lift}(G) \cup {\rm Quad}_B \right> = \left< {\rm Lift}(F') \cup {\rm Lift}(G') \cup {\rm Quad}_B \right>
\end{equation}
where $F'$ and $G'$ are Gr\"obner bases for $I$ and $J$ respectively.  But the ideal on the right hand side of Equation \ref{eqn:gens} is $I \times_\ca J$, since a Gr\"obner basis of an ideal generates that ideal.
\end{proof}

\begin{cor}\label{cor:squarefree}
Suppose that $\ca$ is linearly independent, and that ${\rm in}_{\omega_1} (I)$ and ${\rm in}_{\omega_2}(J)$ are squarefree monomial ideals.  Then ${\rm in}_{\phi^*_B (\omega_1, \omega_2) + \epsilon \omega} (I \times_\ca J)$ is a squarefree monomial ideal.
\end{cor}
\begin{proof}
Since ${\rm in}_{\omega_1}(I)$ and ${\rm in}_{\omega_2}(J)$ are squarefree monomial ideals, then the $\phi^*_B(\omega_1, \omega_2) + \epsilon \omega$ initial term of every polynomial in ${\rm Lift}(F)$ and ${\rm Lift}(G)$ is squarefree.  Since the $\omega$ leading terms of ${\rm Quad}_B$ are also squarefree, and the union of ${\rm Lift}(F)$, ${\rm Lift}(G)$, and ${\rm Quad}_B$ form a Gr\"obner basis for $I \times_\ca J$ we deduce that  ${\rm in}_{\phi^*_B (\omega_1, \omega_2) + \epsilon \omega} (I \times_\ca J)$ is a squarefree monomial ideal.
\end{proof}

It seems natural to ask whether any more refined information about the toric fiber product $I \times_\ca J$ can be computed from  $I$ and $J$ (for instance, homological properties, Betti numbers, Hilbert series).  In view of Example \ref{ex:mono}, it seems that there is no hope of explicitly determining any of these properties from $I$ and $J$,  since the sums of even simple monomial ideals can display complicated behavior.  There is, however, a compact description of the multigraded Hilbert function and series of $I \times_\ca J$.  We refer the reader to \cite{Miller2005} for an introduction to multigraded Hilbert functions.  Note that since we have assumed that the underlying grading is positive, $\dim_\kk ( (\kk[x]/I)_u )$ and $\dim_\kk (( \kk[y]/J)_u)$ are finite for all $u \in \nn \ca$.  Thus the Hilbert function and series exist.

\begin{cor}
Suppose that $\ca$ is linearly independent.  Then the Hilbert functions of $\kk[x]/I$, $\kk[y]/J$ and $\kk[z]/ I \times_\ca J$ satisfy
$$h(\kk[z]/ I \times_\ca J; u) =      h(\kk[x]/I; u) h(\kk[y]/J; u).$$ 
Thus the Hilbert series of $\kk[z]/ I \times_\ca J$ is the Hadamard product
$$H(\kk[z]/ I \times_\ca J; {\bf t}) =      H(\kk[x]/I; {\bf t} ) (\kk[y]/J; {\bf t}).$$
\end{cor}

\begin{proof}
It suffices to produce the first equation, since this is the definition of the Hadamard product of two series.  Let $G$ be a Gr\"obner basis for $I \times_\ca J$ constructed according to Theorem \ref{thm:maingro}.  Let $M = {\rm in}_{ \phi^*_B(\omega_1, \omega_2) + \epsilon \omega}(I \times_\ca J)$ be the corresponding initial ideal.  We need to count the number of monomials in the $u$ graded pieces of $\kk[z] \setminus M$.  A monomial $m = z^\alpha$ belongs to $\kk[z] \setminus M$ if and only if $m$ is not divisible by any of the initial terms of polynomials of ${\rm Lift}(F)$, ${\rm Lift}(G)$ or ${\rm Quad}_B$.  This implies that $\phi_B(m)$ is not divisible by the leading terms of polynomials of $F$ or $G$.  The monomial $\phi_B(m)$ is a product of a degree $u$ monomial in $\kk[x] \setminus {\rm in}_{\omega_1}(I)$ and a degree $u$ monomial in $\kk[y] \setminus {\rm in}_{\omega_2}(J)$.  Any such product has a unique preimage that is not divisible by any leading term in ${\rm Quad}_B$.  Thus we have shown that the standard monomials of $I \times_\ca J$ of degree $u$ are in bijection with a product of a standard monomial of $I$ of degree $u$ and a standard monomial of $J$ of multidegree $u$, which yields the desired equality of Hilbert functions.
\end{proof}

To close this section, we provide an example which shows the necessity of the condition that $\ca$ be a linearly independent in all of the preceding theorems.

\begin{ex} \label{ex:3cycle}
For a fixed vector of positive integers  $d = (d_1, d_2, d_3)$ with $d_i > 1$ for all $i$, consider the ring homomorphism
$$\phi_d :  \kk[p_{i_1i_2i_3} \, \, | \, \, i_j \in [d_j]] \rightarrow \kk[a_{i_1i_2}, b_{i_1i_3} , c_{i_2i_3} \, \, | \, \, i_j \in [d_j] ]$$
$$ p_{i_1i_2i_3} \mapsto a_{i_1i_2} b_{i_1 i_3} c_{i_2 i_3}.$$
Computing the minimal generators of the ideal $K_d = \ker(\phi_d)$ for various values of $d$ is a benchmark problem for algorithms for computing  Gr\"obner bases of toric ideals \cite{Hemmecke2005, Lauritzen2005}.  These ideals are extremely complicated and while there are explicit generating sets known for some special values of $d$, there is, at present, no uniform description of the generating sets of these ideals.  Work of De Loera and Onn \cite{DeLoera2004} suggests that it is impossible for any simple description of a generating set to exist for all $d$.  In particular, it is known that $K_d$ requires minimal generators of degree at least than $2 \min(d_1,d_2,d_3)$

The ideals $K_d$ are examples of toric fiber products that do not have linearly independent $\ca$.  In particular, consider the grading $\deg (p_{i_1i_2 i_3} ) = {\bf e}_{i_1} \oplus {\bf e}_{i_2}$ where ${\bf e}_k$ denotes the $k$-th standard unit vector.  Consider the ring homomorphism
$$\psi_d : \kk[q_{i_1i_2i_3} \, \,  i_j \in [d_j]] \rightarrow \kk[b_{i_1i_3} , c_{i_2i_3} \, \, | \, \, i_j \in [d_j] ]$$
$$q_{i_1i_2i_3} \mapsto  b_{i_1i_3}c_{i_2i_3}$$
and let $I_d = \ker(\psi_d)$.  This ideal is generated by quadrics:  it is an ideal of the form $I_B$ arising from a toric fiber product with linearly independent $\ca$ as in Proposition \ref{prop:quad}.
Grade $\kk[q]$ by $\deg(q_{i_1i_2i_3}) = {\bf e}_{i_1} \oplus {\bf e}_{i_2}$.
Let $J_d = \langle 0 \rangle$ be the zero ideal in $\kk[a]$ with the grading $\deg(a_{i_1i_2}) = {\bf e}_{i_1} \oplus {\bf e}_{i_2}$.  Let $\ca = \{{\bf e}_{i_1} \oplus {\bf e}_{i_2} \, \, | \, \, i_j \in [d_j] \}$.  Then $K_d = I_d \times_\ca J_d$.  So although $I_d$ and $J_d$ are generating in degree $2$ or less, $K_d$ can require minimal generators of arbitrarily large degree.
\end{ex}

%%%%%%%%%%%%%%%%%%%%%%%%%%%%%%%%%%%%%%%
%%%%%%%%%%%%%%%%%%%%%%%%%%%%%%%%%%%%%%%
%%%%%%%%%%%%%%%%%%%%%%%%%%%%%%%%%%%%%%%
%%%%%%%%%%%%%%%%%%%%%%%%%%%%%%%%%%%%%%%

\section{Applications}\label{sec:apps}

\subsection{Segre products}

The simplest example of the toric fiber product is the usual Segre product of two projective schemes.  In this setting $I \subset \kk[x]$ and $J \subset \kk[y]$ are homogeneous ideals in the usual coarse grading by degree.  The monomial homomorphism $\phi_B$ is
$$\phi_B : \kk[z] \rightarrow \kk[x,y]$$
$$z_{jk} \mapsto  x_i y_j.$$

In this case, all variables $x_j$ and $y_k$ have degree $1$ and since $\ca = \{1\}$ is a linearly independent set we are in a position to apply the results of Section \ref{sec:mono}. 
Theorem \ref{thm:main} shows how to produce Gr\"obner bases for $I \times_\ca J$ from Gr\"obner bases for $I$ and $J$.  For instance, we deduce the following corollary, which appears in \cite{Ha2002}.

\begin{cor}
Let $d$ be a vector of nonnegative integers greater than 1.  Let $K_d$ be the defining ideal of the product of projective spaces $\pp^{d_1-1} \times \cdots \times \pp^{d_n-1}$ under the standard Segre embedding into $\pp^{d_1\cdots d_n -1}$.  Then $K_d$ is generated by the $2 \times 2$ minors of all flattenings of a generic $d_1 \times \cdots \times d_n$ tensor.  Furthermore, these $2 \times 2$ minors form a Gr\"obner basis for $K_d$ and the ring $\kk[z]/K_d$ is Cohen-Macaulay.
\end{cor}

A flattening is a matrix obtained by partitioning the $n$ indices of a tensor into two nonempty sets.  For instance, a flattening of a $2 \times 2 \times 2 \times 2$ tensor with partition $\{1,2\}| \{3,4\}$ is
$$
\begin{pmatrix}
z_{1111} & z_{1112} & z_{1121} & z_{1122} \\
z_{1211} & z_{1212} & z_{1221} & z_{1222} \\
z_{2111} & z_{2112} & z_{2121} & z_{2122} \\
z_{2211} & z_{2212} & z_{2221} & z_{2222} \\
\end{pmatrix}.
$$
In \cite{Ha2002}, Ha refers to the $2 \times 2$ minors of flattenings as minors of a box-shaped matrix.

\begin{proof}
It is easy to see that any $2 \times 2$ minor of a flattening belongs to $K_d$ and, furthermore, these are the only degree two binomials in the ideal $K_d$.  Note that $K_d = K_d' \times_\ca \langle 0 \rangle$ where $d' = (d_1,\ldots, d_{n-1})$ and $\ca = \{1\}$.  By repeatedly applying Theorem \ref{thm:main}, we can construct generating sets and Gr\"obner bases for $K_d$.  The lifting operation preserves degrees and every polynomial in ${\rm Quad}_B$ has degree two, so the resulting Gr\"obner bases have degree two.  Since the only degree two binomials in $K_d$ are the $2 \times 2$ minors of flattenings, these must form a Gr\"obner basis.  By Corollary \ref{cor:squarefree}, the resulting initial ideal is squarefree.  Since $K_d$ is a toric ideal, this implies that the $\kk[z]/{\rm in}_\omega(K_d)$ is Cohen-Macaulay (the simplicial complex associated to ${\rm in}_\omega (K_d)$ is a regular triangulation \cite[Theorem 8.3]{Sturmfels1996} and, hence, a shellable ball) which in turn implies that $\kk[z]/ K_d$ is Cohen-Macaulay.
\end{proof}

%%%%%%%%%%%%%%%%%%%%%%%%%%%%%%%%%%%%%%%
%%%%%%%%%%%%%%%%%%%%%%%%%%%%%%%%%%%%%%%
%%%%%%%%%%%%%%%%%%%%%%%%%%%%%%%%%%%%%%%
%%%%%%%%%%%%%%%%%%%%%%%%%%%%%%%%%%%%%%%

\subsection{Reducible models}\label{sec:red}

Probably the first instance where Theorem \ref{thm:main} was used in some generality was in the study of the class of reducible hierarchical models in \cite{Dobra2004} and \cite{Hosten2002}.   Hierarchical log-linear models are a class of statistical models use in the analysis of multivariate discrete data.  To each such hierarchical model is associated a toric ideal $I_{\Delta,d}$.  The generators of the toric ideal $I_{\Delta,d}$ are useful for performing various statistical tests, as first demonstrated in \cite{Diaconis1998}.  In this section, we will only describe these models in a purely algebraic language and show how results about the Gr\"obner bases of reducible models follow from the theory in Section \ref{sec:mono}.

Let $\Delta$ be a simplicial complex with ground set $[n]$ and let $d = (d_1, \ldots, d_n)$ be a vector of integers with $d_i \geq 2$ for all $i$.  We suppose that $| \Delta | = \cup_{F \in \Delta} F = [n]$.  For a subset $F \subset [n]$, we use the notation $D_F$ to denote the set of indices:
$$D_F = \prod_{k \in F}  [d_k].$$ 
For a given string of indices ${\bf i} \in D_{[n]}$, ${\bf i}_F$ is the subvector ${\bf i}_F = (i_{k_1}, \ldots, i_{k_s} )$ where $F =  \{k_1, \ldots, k_s \}$.
Denote by
$$\kk[p] = \kk[p_{\bf i} \,\, | \, \, {\bf i} \in D_{[n]} ]$$
and
$$\kk[a] = \kk[ a^F_{{\bf j}_F} \, \, | \, \, F =  \in {\rm facet}(\Delta), {\bf j}_F \in D_F ].$$ 
Consider the ring homomorphism
$$\phi_{\Delta,d} :  \kk[p] \rightarrow \kk[a]$$
$$p_{\bf i} \mapsto  \prod_{F  \in {\rm facet}(\Delta) }  a^F_{{\bf i}_F}.$$

\begin{df}
The toric ideal $I_{\Delta,d} = \ker( \phi_{\Delta,d})$ is the ideal of the hierarchical model defined by $\Delta$ and $d$.
\end{df}

\begin{ex}\label{ex:red}
Let $\Delta$ be the simplicial complex with facets ${\rm facet}(\Delta) = \{\{1,2\},\{2,3\}, \{3,4\}\}$.  Then
$$\phi_{\Delta,d} :  \kk[p_{i_1i_2i_3i_4} : i_l \in [d_l] ] \rightarrow \kk[a^{\{1,2\}}_{j_1j_2}, a^{\{2,3\}}_{j_2j_3}, a^{\{3,4\}}_{j_3j_4} \, \, | \, \, j_l \in [d_l] ]$$
$$p_{i_1i_2i_3i_4} \, \, \,  \longmapsto \, \, \,   a^{\{1,2\}}_{i_1i_2} a^{\{2,3\}}_{i_2i_3} a^{\{3,4\}}_{i_3i_4}. $$
If $\Delta$ is the simplicial complex with facets ${\rm facet}(\Delta) = \{\{1,2\},\{1,3\},\{2,3\}\}$, then we get the ring homomorphism from Example \ref{ex:3cycle}. \qed
\end{ex}

\begin{df}
A simplical complex $\Delta$ is called \emph{reducible} if there are two subcomplexes $\Delta_1, \Delta_2 \subset \Delta$ such that $\Delta_1 \cup \Delta_2 = \Delta$ and $\Delta_1 \cap \Delta_2 = 2^S$ for some $S \subset [n]$.  The set $S$ is called a \emph{separator}.
\end{df}

For instance, the first simplicial complex from Example \ref{ex:red} is reducible with $S = \{1\}$ or $S = \{2\}$, whereas the second simplicial complex from Example \ref{ex:red} is not reducible.   We will show that if $\Delta$ is reducible, the ideal $I_{\Delta,d}$ can be written as a toric fiber product.  This will allow us to deduce Theorem 4.17 from \cite{Hosten2002}.  

To this end, let $d^1$ and $d^2$ be the induced subvectors of $d$ on the index sets $| \Delta_1 |$ and $| \Delta_2 |$ respectively.  That is $d^1 = d_{|\Delta_1|}$ and $d^2 = d_{|\Delta_2|}$.
For $l \in \{1,2\}$ let 
$$\kk[p]_l  = \kk [ p_{{\bf i}_{|\Delta_l |}} \, \, | \, \, {\bf i }_{|\Delta_l |} \in  D_{|\Delta_l |} ]$$
and consider the ring homomorphism 
$$\phi_{\Delta_l, d^l } : \kk[p]_l \rightarrow \kk[a]$$
$$ p_{{\bf i}_{|\Delta_l |} } \mapsto  \prod_{F \in {\rm facet}(\Delta_l)}  a^F_{({\bf i}_{|\Delta_l |} )_F}. $$

We denote by $I_{\Delta_l,d^l}$ the kernel of $\phi_{\Delta_l, d^l}$, which is the toric ideal of the hierarchical model associated to $\Delta_l$.  Since $\Delta_1 \cap \Delta_2 = 2^S$, there exists two facets, $F_1 \in \Delta_1$ and $F_2 \in \Delta_2$ such that $S \subseteq F_1$ and $S \subseteq F_2$.  We introduce a grading on $\kk[a]$ so that 
$$\deg(a^{F}_{{\bf j}_{F}}) = \left\{ \begin{array}{cl}
e_{({\bf j}_{F})_S} & \mbox{ if } F \in \{F_1, F_2\} \mbox{ and } \\
 0 &  \mbox{ otherwise } \end{array} \right.$$
 The vector $e_{{\bf i}_S}$ for ${\bf i}_S \in d_S$ is the standard unit vector in $\zz^{D_s}$ with a $1$ in the ${\bf i}_S$ position and a zero elsewhere.  This multigrading on $\kk[a]$ induces a multigrading on $\kk[p]$, $\kk[p]_1$ and $\kk[p]_2$, where, for instance, we take the degree of $p_{\bf i}$ to be $\deg(p_{\bf i}) = e_{{\bf i}_S}$.  Thus, all of $I_{\Delta,d}$, $I_{\Delta_1, d^1}$ and $I_{\Delta_2, d^2}$ are homogeneous with respect to this multigrading, because the maps $\phi_{\Delta,d}$, $\phi_{\Delta_1,d^1}$ and $\phi_{\Delta_2, d^2}$ all preserve the multidegree.

\begin{thm}\label{thm:red}
Let $\Delta$ be reducible with components $\Delta_1$ and $\Delta_2$ and separator $S$.  Then
$$I_{\Delta, d} = I_{\Delta_1, d^1} \times_\ca I_{\Delta_2, d^2}$$
with $\ca$ linearly independent.
\end{thm}

\begin{proof}
It is clear that $\ca$ is linearly independent since it is a collection of disjoint standard unit vectors.  Suppose that $T$ is an arbitrary face of $\Delta$.
We first note the general fact that if we modify the ring homomorphism $\phi_{\Delta, d}$ so that
$$p_{\bf i} \mapsto    a^T_{{\bf i}_T}  \cdot \prod_{F \in { \rm facet}(\Delta)} a^F_{{\bf i}_F}$$
this does not change $\ker(\phi_{\Delta,d})$, since any variable $a^T_{{\bf i}_T}$ must appear with precisely same multiplicitly as $a^F_{{\bf i}_F}$ for any facet $F$ with $T \subseteq F$ in the image of a monomial $\phi_{\Delta,d}$.   Consider the modified parametrization
$$p_{\bf i} \mapsto    (a^S_{{\bf i}_S})^2  \cdot \prod_{F \in { \rm facet}(\Delta)} a^F_{{\bf i}_F}$$
where $S$ is the separator.  We factorize the expression on the right as:
$$(a^S_{{\bf i}_S})^2  \cdot \prod_{F \in { \rm facet}(\Delta)} a^F_{{\bf i}_F}  \quad = \quad 
\prod_{l \in \{1,2\} }  \left( a^S_{{\bf i}_S}  \cdot  \prod_{F \in {\rm facet} (\Delta) \cap {\rm facet}(\Delta_l)} a^F_{{\bf i}_F} \right) .
$$
Since the expression on the right inside the parentheses involves a product of terms for all facets of $\Delta_l$, this product represents the parametrization for $I_{\Delta_l, d^l}$.  In other words, we have shown that the ring homomorphism $\phi_{\Delta,d}$ factors through
$$\phi_B:  \kk[p] \rightarrow \kk[p]_1 \otimes \kk[p]_2$$
$$p_{\bf i} \mapsto  p_{{\bf i}_{|\Delta_1|}} \otimes p_{{\bf i}_{|\Delta_2|}}.$$
Thus $I_{\Delta, d}$ is a toric fiber product.
\end{proof}

%%%%%%%%%%%%%%%%%%%%%%%%%%%%%%%%%%%%%%%
%%%%%%%%%%%%%%%%%%%%%%%%%%%%%%%%%%%%%%%
%%%%%%%%%%%%%%%%%%%%%%%%%%%%%%%%%%%%%%%
%%%%%%%%%%%%%%%%%%%%%%%%%%%%%%%%%%%%%%%

\subsection{Reducible models with hidden variables}\label{sec:redhidden}

Pushing the idea from Section \ref{sec:red} one step further, we can also use the machinery to compute the ideals of reducible models from submodels when some of the random variables are hidden.  Parametrically, we have the following setup for the algebraic description of a hidden variable models.  

Let $\Delta$ be a simplicial complex on $[n]$, $d = (d_1, \ldots, d_n)$ and vector of integers with $d_i \geq 2$ for all $i$.  Let $H = \{h_1, \ldots, h_t \}  \subset [n]$ be the collection of \emph{hidden nodes} and $O = [n] \setminus H$ be the collection of \emph{observed nodes}.  
Let
$$\kk[q] = \kk[q_{{\bf i}} \, \, | \, \, {\bf i}_O \in D_O, i_l = \bullet \mbox{ if } l \in H ]$$
and let $\kk[p]$, $\phi_{\Delta,d}$, and $I_{\Delta,d}$ be defined as in Section \ref{sec:red}.
If ${\bf i}_O \in D_O$ and ${\bf j}_H \in D_H$ we use the notation $p_{{\bf i}_O {\bf j}_H}$ to denote the indeterminate $p_{{\bf i}}$ such that ${\bf i}_l = ({\bf i}_O)_l$ if $l \in O$ and ${\bf i}_l = ({\bf j}_H)_l$ if $l \in H$.  Consider the ring homomorphism 
$$\psi_H :  \kk[q] \rightarrow \kk[p]$$
$$q_{\bf i} \mapsto  \sum_{{\bf j}_H \in D_H}   p_{{\bf i}_O {\bf j}_H }.$$

Denote by $I_{H,\Delta,d} =  \ker ( \psi_H \circ \phi_{\Delta,d} ) $ which is the ideal of the hidden variable hierarchical model.   This ideal is rarely a toric ideal.  Hidden variable graphical models have been studied in \cite{Garcia2005} and \cite{Garcia2004} from the perspective of computational algebra though we seem to be the first to write down some general principles for determining their defining prime ideals. 

\begin{df}
We call the hidden variable ideal $I_{H,\Delta, d}$ \emph{reducible} if $\Delta$ is a reducible simplicial complex and $H \cap S = \emptyset$.
\end{df}

Suppose that $I_{H,\Delta,d}$ is a reducible hidden variable ideal.  Let $\Delta_1$ and $\Delta_2$ be the two component subcomplexes, let $H_1 = H \cap |\Delta_1|$ and $H_2 = H \cap |\Delta_2|$ and let $d^1$ and $d^2$ be the two induced vectors of indices.  Denote by $\kk[q]_1$ and $\kk[q]_2$ the two rings with variables indexed by the elements of $D_{|\Delta_1|}$ and $D_{|\Delta_2|}$.

\begin{thm}\label{thm:hidred}
Let $I_{H,\Delta,d}$ be a reducible hidden variable ideal. Let $S$ be the separator.  In each of the rings $\kk[q]$, $\kk[q]_1$, and $\kk[q]_2$, let the degree of a variable be $\deg(q_{\bf i}) = e_{{\bf i}_S}$.  Then
$$I_{H,\Delta,d} = I_{H_1, \Delta_1, d^1} \times_\ca I_{H_2, \Delta_2, d^2}$$
with $\ca$ linearly independent.
\end{thm}

\begin{proof}
The same proof as of Theorem \ref{thm:red} applies here.
\end{proof}

As an example of an application of Theorem \ref{thm:hidred} we will deduce a Corollary generalizing Example \ref{ex:dets}.  We call this example the \emph{partially hidden Markov chain}.
Let $\Delta = \{\{1,2\},\{2,3\}, \{3,4\}, \ldots, \{2n,2n+1\} \}$ be a chain of odd length and suppose that $H = \{2,4,6, \ldots, 2n\}$ consists of all the even numbers.  To describe the generators of the ideal $I_{H,\Delta,d}$ we need two matrix constructions.  First, each even number $2j \in [2n+1]$ defines a flattening of the $n+1$ tensor $(q_{\bf i})$ into a matrix $X_j$.  The rows and column indices of the matrix $X_j$ are the elements of $D_F$ with $F = \{1,3,\ldots, 2j-1\}$ and the elements of $D_{F'}$ with $F' = \{2j+1, 2j+3, \ldots, 2n+1 \}$, respectively.  Thus $X_j$ is a $d_1d_3 \cdots d_{2j-1} \times  d_{2j+1}d_{2j+3} \cdots d_{2n+1}$ matrix.  The entry in the ${\bf i}_F$ row and ${\bf j}_{F'}$ column is $q_{{\bf i}_F {\bf j}_{F'}}$ (with appropriate $\bullet$'s added).

Second, to each odd number $2j+1 \in [2n+1]$ with $0< j < n$, and each $i \in [d_{2j+1}]$  we introduce a matrix $Y_{j,i}$ which is a flattening of an $n$ dimensional slice of the $(n+1)$-dimensional tensor $(q_{\bf i})$.  The row and column indices of the matrix $Y_{j,i}$ are the elements of $D_F$ with $F = \{ 1,3, \ldots, 2j-1\}$ and the elements of $D_{F'}$ with $F' = \{2j+3, 2j+5, \ldots, 2n+1\}$.  
Thus, $Y_{j,i}$ is a $d_1d_3\cdots d_{2j-1} \times d_{2j+3} d_{2j+5} \cdots d_{2n+1}$ matrix.  The entry in the ${\bf i}_F$ row and ${\bf j}_{F'}$ column is $q_{{\bf i}_F  i {\bf j}_{F'}}$ (with appropriate bullets added.  Examples of $X_j$ and $Y_{j,i}$ are illustrated in Example \ref{ex:dets} (note that the second $3 \times 9$ matrix in Example \ref{ex:dets} is the transpose of $X_4$).

\begin{cor}\label{cor:chain}
Let $\Delta = \{\{1,2\},\{2,3\}, \{3,4\}, \ldots, \{2n,2n+1\} \}$ be a chain of odd length and suppose that $H = \{2,4,6, \ldots, 2n\}$ consists of all the even numbers. Let $G$ be the union of all $(d_{2j} +1) \times (d_{2j} + 1)$ minors of $X_j$ for $j \in [n]$ such that $d_{2j} < \min(d_{2j-1}, d_{2j+1})$ together with the union of all the $2 \times 2 $ minors of $Y_{j,i}$ for $0 < j < n$ and $i \in [d_{2j+1}]$.  Then $G$ is a Gr\"obner basis for $I_{H,\Delta, d}$.
\end{cor}

\begin{proof}
Note that $I_{H,\Delta, d}$ is reducible with separator $S = \{2n\}$.  Thus, by applying induction and liberal application of Theorem \ref{thm:hidred}, we need only determine Gr\"obner bases for $I_{H',\Delta', d'}$ in the special case of $\Delta' = \{\{1,2\},\{2,3\}\}$ and $H' = {2}$ for all triples $d' = (d_1, d_2, d_3)$.  Then we can lift the polynomials to get a Gr\"obner basis for $I_{\Delta, H, d}$.  
In this special case, we are considering the ideal that is the kernel of the ring homomorphism
$$\phi_{d} : \kk[q_{i_1 \bullet i_3} \, | \, i_1 \in [d_1], i_3 \in [d_3] ] \rightarrow  \kk[a_{i_1i_2}, b_{i_2i_3} \, | \, i_j \in [d_j]  ]$$
$$ q_{i_1 \bullet i_3} \mapsto   \sum_{i_2 = 1}^{d_2}  a_{i_1 i_2} b_{i_2 i_3}.$$
 
Thus, $I_{H',\Delta',d}$ is the vanishing ideal of $Sec^{d_2-1}( \pp^{d_1 -1} \times \pp^{d_3 - 1})$, the variety of secant $(d_2 - 1)$-planes to the Segre embedding of $\pp^{d_1 -1} \times \pp^{d_3 -1}$.  If $d_2 \geq \min(d_1,d_3)$ then $I_{H,\Delta,d} =\left< 0 \right>$ is the zero ideal.  If $d_2 < \min(d_1,d_3)$  then $I_{H,\Delta,d}$ is generated by the $(d_2 +1) \times (d_2 +1)$ minors of the matrix $(q_{i_1 \bullet i_3})$.  Lifting these determinantal polynomials to $\kk[q]$ yields the minors of 
matrix $X_j$.   The minors of the matrices $Y_{j,i}$ are the elements of ${\rm Quad}_B$ for each of the toric fiber products that are used in building up $I_{H, \Delta, d}$.
\end{proof}

%%%%%%%%%%%%%%%%%%%%%%%%%%%%%%%%%%%%%%%
%%%%%%%%%%%%%%%%%%%%%%%%%%%%%%%%%%%%%%%
%%%%%%%%%%%%%%%%%%%%%%%%%%%%%%%%%%%%%%%
%%%%%%%%%%%%%%%%%%%%%%%%%%%%%%%%%%%%%%%

\subsection{Group-based models on phylogenetic trees} 

In this section, we show how the toric fiber product arises in the construction of phylogenetic invariants for group-based models on phylogenetic trees.  The fact that these phylogenetic models are toric fiber products plays a significant role in \cite{Sturmfels2005} where phylogenetic invariants for the group based models were originally constructed.  The toric fiber product also plays a prominent role in \cite{Buczynska2006} where, in the case of trivalent trees with underlying group $\zz_2$, it is proven that these phylogenetic ideals are Gorenstein, their Hilbert polynomials are computed, and their deformations are studied.  For the sake of simplicity of exposition, we will describe the underlying models in the Fourier coordinates (as opposed to the probability coordinates) and we only address the case of models whose labeling function is the identity map.  We refer the reader to \cite{Sturmfels2005} and \cite{Szekely1993} for descriptions of these models in the probability coordinates and the application of the discrete Fourier transform.  Also \cite{Sturmfels2005} contains the full description of these models with arbitrary friendly labeling functions.

Let $T$ be a tree with $n+1$ leaves.  Label the leaves $1, 2, \ldots, n+1$, let the root of $T$ be at the leaf $n+1$, and direct the edges of $T$ away from the root.  Given an edge $e$, a leaf $l$ is called an descendant of $e$ if there is a directed path from $e$ to $l$.  Denote by ${\rm de}(e)$ the set of all descendants of the edge $e$.  We assume that the tree $T$ has the property that for every edge $e$ the set of descendants ${\rm de}(e)$ is an interval of integers, i.e. ${\rm de}(e) = \{i,i+1, \ldots j-1, j\}$ for some $i \leq j \in [n]$.  This amounts to saying that $T$ has a drawing in the plane so that the leaves of $T$ lie on a circle in numerical order.

Let $G$ be a group.  We will use additive notation for $G$ although $G$ might not be abelian. If $g_1, \ldots g_n$ is a sequence of elements in $G$, we denote by $g_e$ the sum of the group elements $g_i$ such that $i \in {\rm de}(e)$; that is,
$$g_e = \sum_{i \in {\rm de}(e)} g_i.$$
Since $G$ may not be abelian we use the convention that the sum is always taken in increasing order of the indices $i \in {\rm de}(e)$.  Let 
$$\kk[q] = \kk[q_{g_1 \ldots g_n } \, \, | \,\, g_i \in G ] \mbox{ and }  \kk[a] = \kk[a^{(e)}_{h} \, \, | e \in E(T), h \in G]$$
and consider the ring homomorphism
$$\phi_{G,T} : \kk[q]  \rightarrow \kk[a]$$
$$ q_{g_1 \ldots g_n} \mapsto  \prod_{e \in E(T)}   a^{(e)}_{g_e}.$$

\begin{df}
The ideal $I_{G,T} = \ker(\phi_{G,T})$ is the ideal of the group-based phylogenetic model with group $G$ and tree $T$.
\end{df}

Note that since $T$ is an acyclic directed graph, there is an induced partial order on the edges of $T$.  Namely $e < e'$ if there is a directed path from $e'$ to $e$.   Let $T$ be a tree that contains an interior edge (an edge not incident to any leaf).  Then $e$ induces a decomposition of $T$ as the composition $T_e^+ \ast T_e^-$ where $T_e^-$ is the subtree of $T$ consisting of all edges $e' \in T$ with $e' \leq e$ and $T_e^+$ consists of all edges $e' \in T$ such that $e' \not < e$.  Thus $T_e^-$ and $T_e^+$ overlap in the single edge $e$.  We root $T_e^-$ by the tail of $e$, and keep the root of $T_e^+$ at the original root $n+1$.

Without loss of generality, we may assume that the nonroot leaves of $T_e^-$ consist of $\{1,2, \ldots, k\}$, and the nonroot leaves of $T_e^+$ are $\{e, k+1, \ldots, n\}$.  Let $\kk[q]_+$ and $\kk[q]_-$ denote the ambient polynomial rings of $I_{G,T^+_e}$ and $I_{G, T^-_e}$, respectively.  

\begin{thm}
Let $T$ be a tree with an interior edge $e$, and resulting decomposition $T = T^+_e \ast T^-_e$.  For each variable $q_{\bf g}$ in $\kk[q]$, $\kk[q]_+$ and $\kk[q]_-$, let $\deg(q_{\bf g})  = e_{g_e}$, the standard unit vector with label $g_e$.  
Then
$$I_{G,T} = I_{G,T^+_e}  \times_\ca I_{G,T^-_e}$$
with $\ca$ linear independent.
\end{thm}

\begin{proof}
Clearly $\ca$ is linearly independent since it consists of standard unit vectors.  To prove that $I_{G,T}$ is a toric fiber product, we use our standard technique of modifying the parameterization.  As usual, it does not hurt to square a variable $a^{(e)}_{g_e}$ everywhere it appears.  Thus we have:
\begin{eqnarray*}
\phi_{G,T} (q_{\bf g}) &  =  &  a^{(e)}_{g_e}  \cdot  \prod_{e' \in E(T)}   a^{(e')}_{g_{e'}} \\
  & = &   \prod_{e' \in E(T^+_e)}  a^{(e')}_{q_{e'}}   \prod_{e' \in E(T^-_e)}  a^{(e')}_{q_{e'}} \\
  & = &  \phi_{G,T^+_e}( q_{{\bf g}_+}) \phi_{G,T^-_e}( q_{{\bf g}_-})
\end{eqnarray*}
Thus, $I_{G, T}$ is a toric fiber product.
\end{proof}

This allows us to deduce the main result from \cite{Sturmfels2005}.

\begin{cor}
If $T$ is a trivalent tree the generators of $I_{G,T}$ can be explicitly determined from the generators of $I_{G, K_{1,3}}$ where $K_{1,3}$ is the three leaf claw tree. 
\end{cor}
\begin{proof}
If $T$ is a trivalent tree, it can be successively be decomposed by the $\ast$ operation until each component tree is a $K_{1,3}$.
\end{proof}

%%%%%%%%%%%%%%%%%%%%%%%%%%%%%%%%%%%%%%%
%%%%%%%%%%%%%%%%%%%%%%%%%%%%%%%%%%%%%%%
%%%%%%%%%%%%%%%%%%%%%%%%%%%%%%%%%%%%%%%
%%%%%%%%%%%%%%%%%%%%%%%%%%%%%%%%%%%%%%%

\end{document}